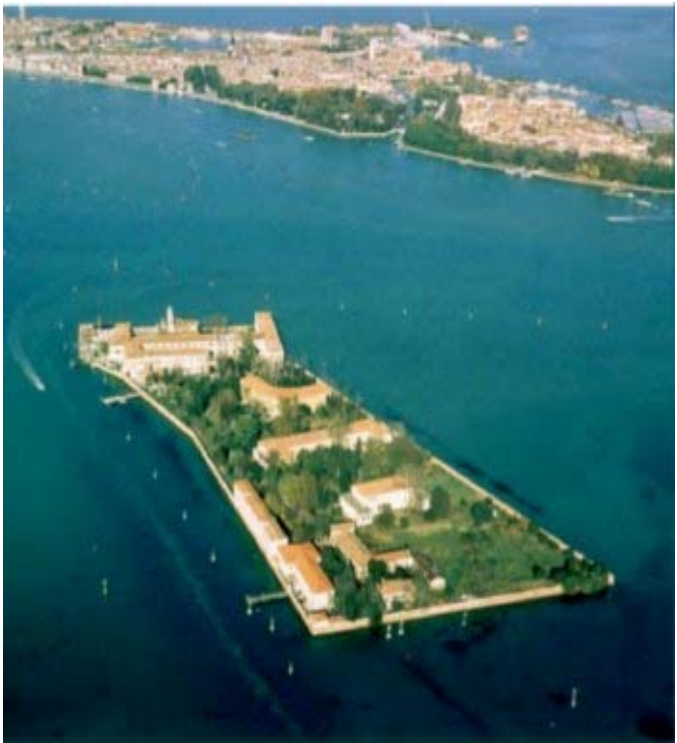

**SIMAI 2004**

modelli e metodi matematici applicati all'industria, alla tecnologia, alla biologia, alla finanza, all'ambiente

Venezia, Isola di San Servolo
20-24 settembre 2004

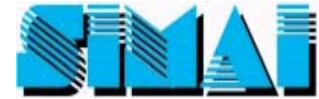

## Using sparse matrices and splines-based interpolation in *computational fluid dynamics* simulations

**Gianluca Argentini**
*Laboratorio di Computazione Avanzata*
Riello Group, Legnago (VR), Italy

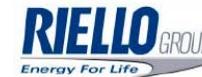

gianluca.argentini@riellogroup.com

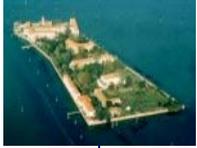

# Summary

- *Position of the computational problem*
- *The method of splines-based interpolations*
- *The computation of the splines coefficients*
- *The computation of the splines values*
- *A principle of virtual equivalence*
- *Pros and cons of the method*

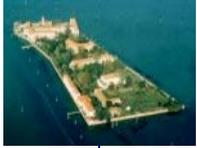

# Numerical problems for burners

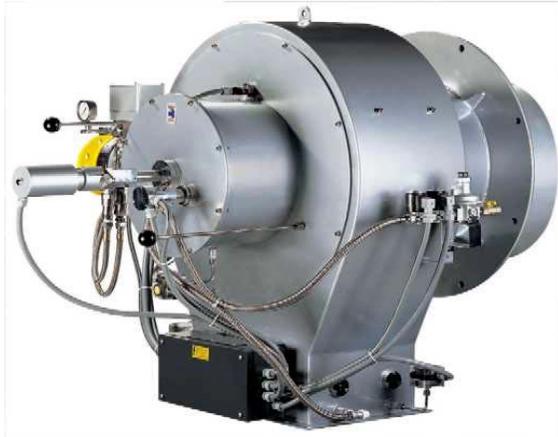

Design, development and engineering of industrial power burners have strong mathematical requests:

• computation of fields of temperature, pressure and velocity in the *combustion chamber*

• correct design of the *combustion head* for an optimal efficiency of the flame

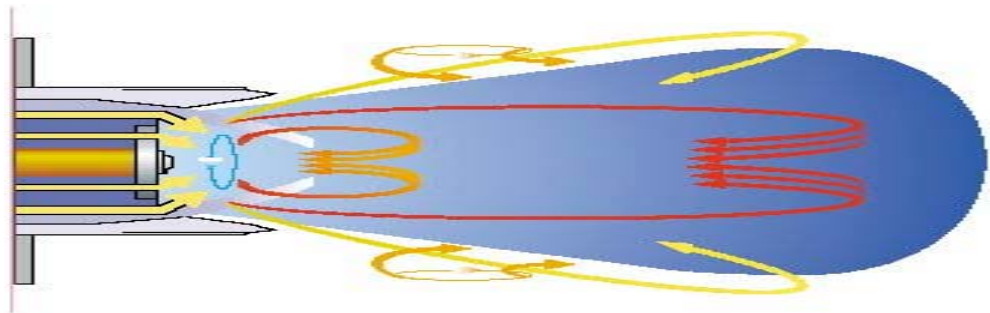

• computation of all the flows (*air, oil, residual gases*) in the burner components

• design of the optimal shape for *ventilation fans*

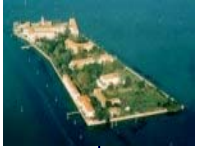

# Limited area meteo model

Two recent floodings by violent meteorological events have induced to consider a service of weather forecasting

*vaghissimi e famosi colli Euganei*

**F. Petrarca, 1304-1374,
700th anniversary of the birth**

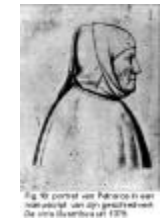

the company has factories in two opposite side of the geographical zone interested by Euganei Hills

Numerical Weather Prediction: based on mathematical models and numerical resolution of a PDEs system on the atmospheric fluid

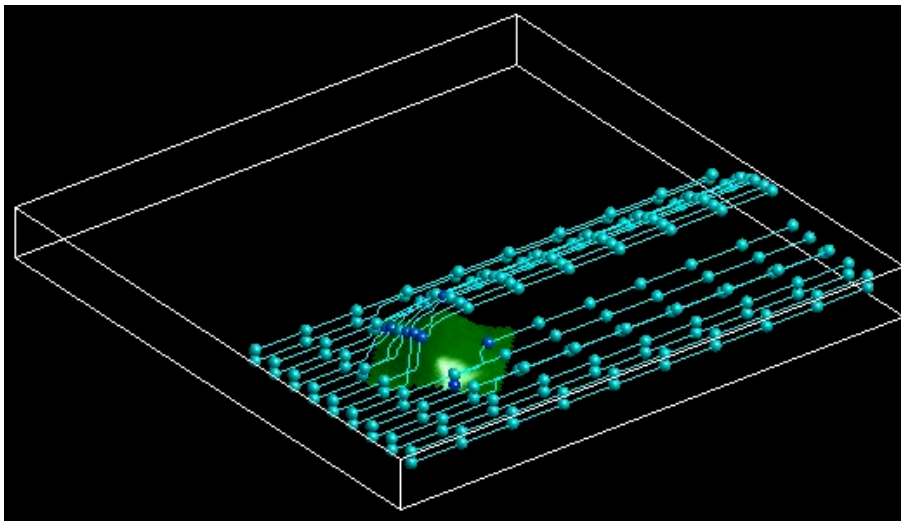

*Streamlines over a hill from a raw version of EHLAM, Euganean Hills Limited Area Model*

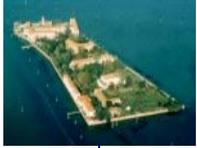

# Computational fluid dynamics

*Numerical resolution* of a system of PDEs for fluid flows is required:

- Navier-Stokes for velocity and pressure of flows

- diffusion-like equation for temperature field

- conservation law for mass and energy of multiphase fluids, i.e liquid-gas oil components

- boundary conditions for the geometry of domain, e.g. combustion chamber

The computational model for burners is quite similar to that of *NWP*  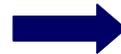  Large use of distributed and parallelized computations on *multiprocessor computers*

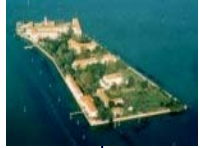

# *Computational complexity analysis for a flow*

***Simple example*** for a detailed knowledge of the velocity-field of fluid particles in the combustion chamber, using a *flow-like grid* :

• **M** is the number of flow streamlines to compute

• **S** is the number of geometrical points for every streamline

High values for **M** are important for a *realistic simulation* of the flow, high values for **S** are important for a fine *graphic resolution* : good values are of order $10^3$ - $10^4$

Using *finite difference* in a computational grid, for every time step the number of computational flops is of order $10^9$, and from Courant-Friedrichs-Lewy (CFL) condition, the time step must to be very small:

*computation and graphic rendering of one minute of flow is very CPU expensive (some Gflop/s) and RAM consuming (hundreds of Mbytes)*

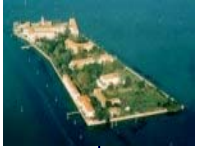

## *Mathematical models and software*

We have experimented three ways:

1. *commercial software*, based on finite elements method for a numerical resolution of Navier-Stokes equations; in general, the accuracy of solution is good, but the methods are not easy customizable and CPU-RAM expensive;

2. *cellular automata model* for the computation of velocity field, based on C or Fortran programs, very useful for generic geometries but RAM consuming; the treatment of the temperature is difficult;

3. *finite difference schema* for the complete system of equations, based on MATLAB or Fortran programs; the computation of the flow is complete, but for a realistic simulation we must verify CFL condition and stability criteria.

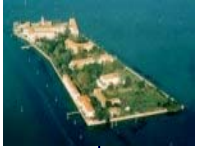
# Cellular automaton model

S. Wolfram (1986) has shown the equivalence between some cellular automata for fluids and Navier-Stokes equations

With some simple customized models of cellular automata we have obtained good geometrical description of flows, but we have noticed difficulties on:

1. correct treatment of boundary conditions;

2. computation of the temperature field;

3. huge consumption of CPUs and RAM.

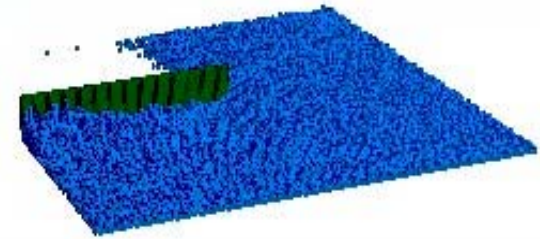

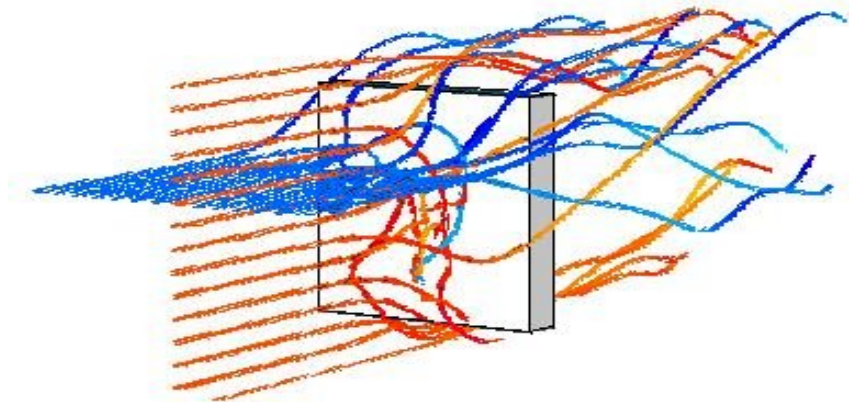

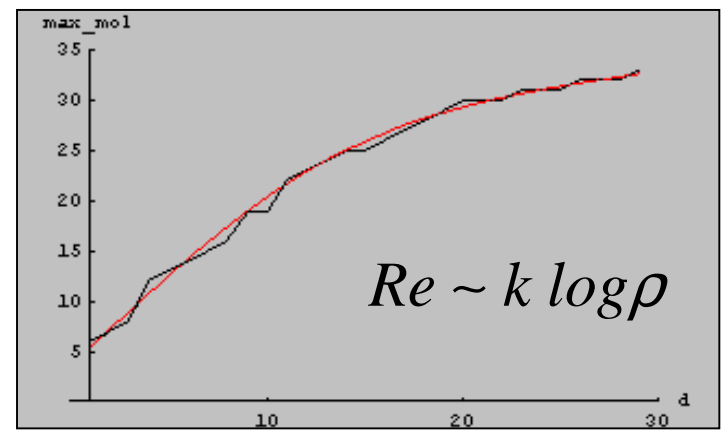

$$Re \sim k \, log\rho$$

# *Finite difference schema model*

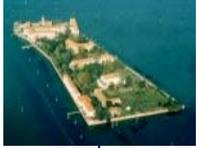

PDEs system from Navier-Stokes equations, mass conservation law, first thermodynamics principle, fluid equation of state

unknowns: velocity vector, pressure, density, temperature

a generic shape of a burner component or the orography of a hill require a non-uniform computational grid

FD schema: Lax-Friedrichs (forward in time, centered in space)

*Example of computational complexity for a **small** simulation.* 3D box 50x25x25 cm, medium scalar velocity of fluid in each cartesian direction of combustion chamber 50 cm/sec, a space resolution of 0.5 cm: *what is a right time-step?* From CFL condition we have

       time-step < 0.5 cm / 50 cm/sec = 0.01 sec

➡ *$O(10^{10})$ flops and **5 GB RAM** for **1** real minute of simulation*

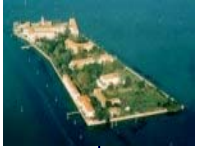

# The method of FD and interpolations

Two problems:

1. time-step is too small and generates a lot of non useful snapshots per second;

2. RAM occupation is very large even in the case of limited simulations.

Suppose to accept 10 snapshots/second; from CFL condition we have

*min space-step* = 0.1 sec * 50 cm/sec = 5 cm

For realistic resolution of single components and good graphic rendering, this value can be too high: for better final results, we have developed a method based on the *interpolation* of the computed values of FD solutions; we have experimented that CPUs effort and RAM occupation are lower than in the case of a fine grid simulation, without significant loss in the final resolution.

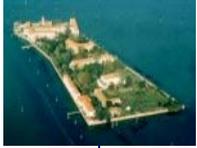
# Phases for interpolations

The method is based on two steps for the graphic rendering of fluid particles trajectories, after the numerical computation of Navier-Stokes or CA model:

1. interpolation by cubic splines of the geometric positions of the particles:

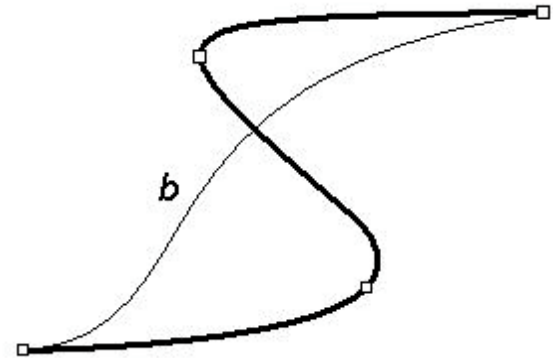

2. fine valuation of every cubics in a suitable set of time values $t_i$:

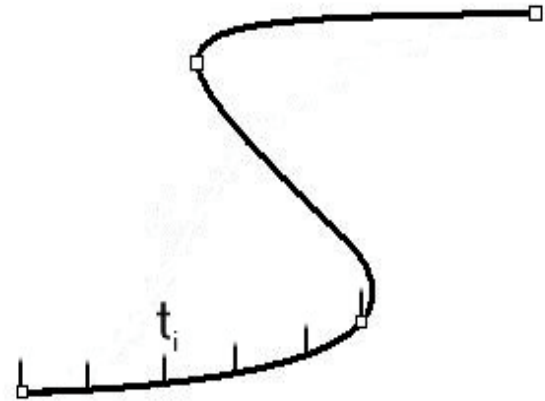

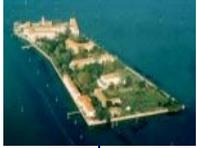

## *Fitting the trajectories*

Let **S** the number of computed velocity vectors in a particle trajectory, **M** the number of trajectories.

What method for interpolation of speed-points?

• Bezier-like is not realistic for rendering the divergence of velocity field

• Chebychev or Least-Squares-like are too rigid in the case of a customized application

• polinomial-like is simple but often shows spurious effects as Runge phenomenon, p.e. :

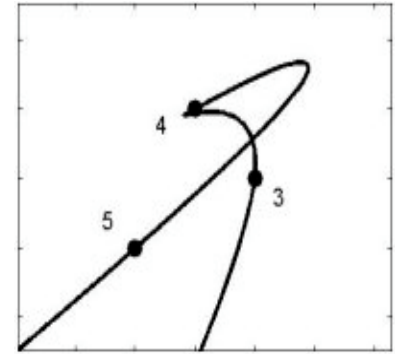

We have obtained better results with a particular splines-based method.

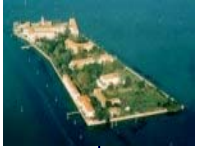

# The splines-based algorithm

Let **S** = 3 x **N** : a trajectory is divided into **N** groups, each of **4** points

At every group the points are interpolated by three cubic polynomials imposing four analytical conditions:

- passage at $P_k$ point, $1 \leq k \leq 3$
- passage at $P_{k+1}$ point
- continuous slope at $P_k$ point
- continuous curvature at $P_k$ point

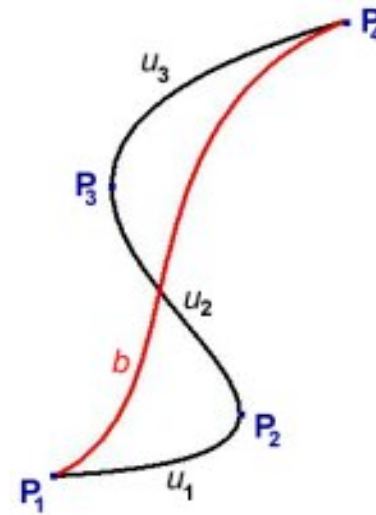

For smooth rendering and for avoiding excessive twisting of trajectories, the cubics $u_k$ are added to the Bezier curve $b$ associated to the four points:

$$v = \alpha b + \beta u_k \qquad 0 < \alpha, \beta < 1$$

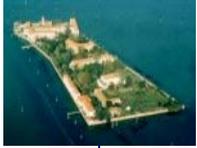

# Finding the splines

We consider $\alpha = \beta = 0.5$

Let $b = As^3 + Bs^2 + Cs + D$ $(0 \leq s \leq 1)$ the Bezier curve of control points $P_1,\ldots,P_4$; for every spline

$u_k = at^3 + bt^2 + ct + d$ $(0 \leq t \leq 1)$

the coefficients can be computed by the system

$$(a, b, c, d) = \mathbf{T}\,(P_{k+1}, P_k, B, C, 1)$$

(matrix-vector multiplication) where the matrix $\mathbf{T}$ is *constant* :

$$\mathbf{T} = \begin{pmatrix} 1 & -1 & -3 & -1 & -6 \\ 0 & 0 & 1 & 0 & 3 \\ 0 & 0 & 2 & 1 & 3 \\ 0 & 1 & 0 & 0 & 0 \end{pmatrix}$$

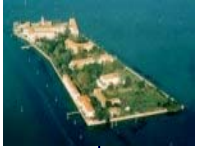

# A global matrix for splines

We define the **G** matrix:

$$G = \begin{pmatrix} T & 0 & . & . & 0 \\ 0 & T & . & . & 0 \\ . & . & . & . & . \\ . & . & . & . & . \\ 0 & 0 & . & . & T \end{pmatrix}$$

(**0** is the 4 x 5 zero-matrix)

**G** is a 4**M** x 5**M** *sparse* matrix with density number < 1/**M**

if **b** = ($P_{k+1}$, $P_k$, $B_1$, $C_1$, 1, . . ., $P_{k+1}$, $P_k$, $B_M$, $C_M$, 1)
we can compute for every two-points group the coefficients of cubic splines for all the **M** trajectories:

**coeff** = **G b**

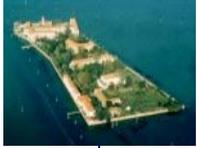

# *Flops and time execution for the splines*

The theoric number of floating point operations for computing the coefficients of all the splines in flow is $O(10\, M^2\, N)$

in a useful *not too small* simulation the **N**, **M** values are of order $10^3$ -$10^4$

the total number of flops is $O(10^{12})$

With a processor having a clock frequency of GHz order the total time can require *some hundreds of seconds*, which is a performance not very good for fast graphics; using

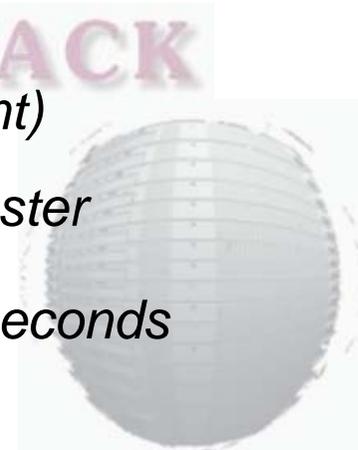

- *some mathematical libraries as LAPACK routines (Fortran calls or Matlab environment)*
- *distributed computation on a multinode cluster*

we have reached a computation time of *some tens of seconds*

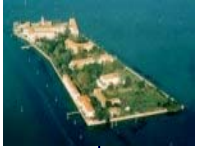

# Computing with Lapack

*Example* : Matlab has internal Lapack level 3 BLAS routines for fast matrix-matrix multiplication and treatment of sparse matrices

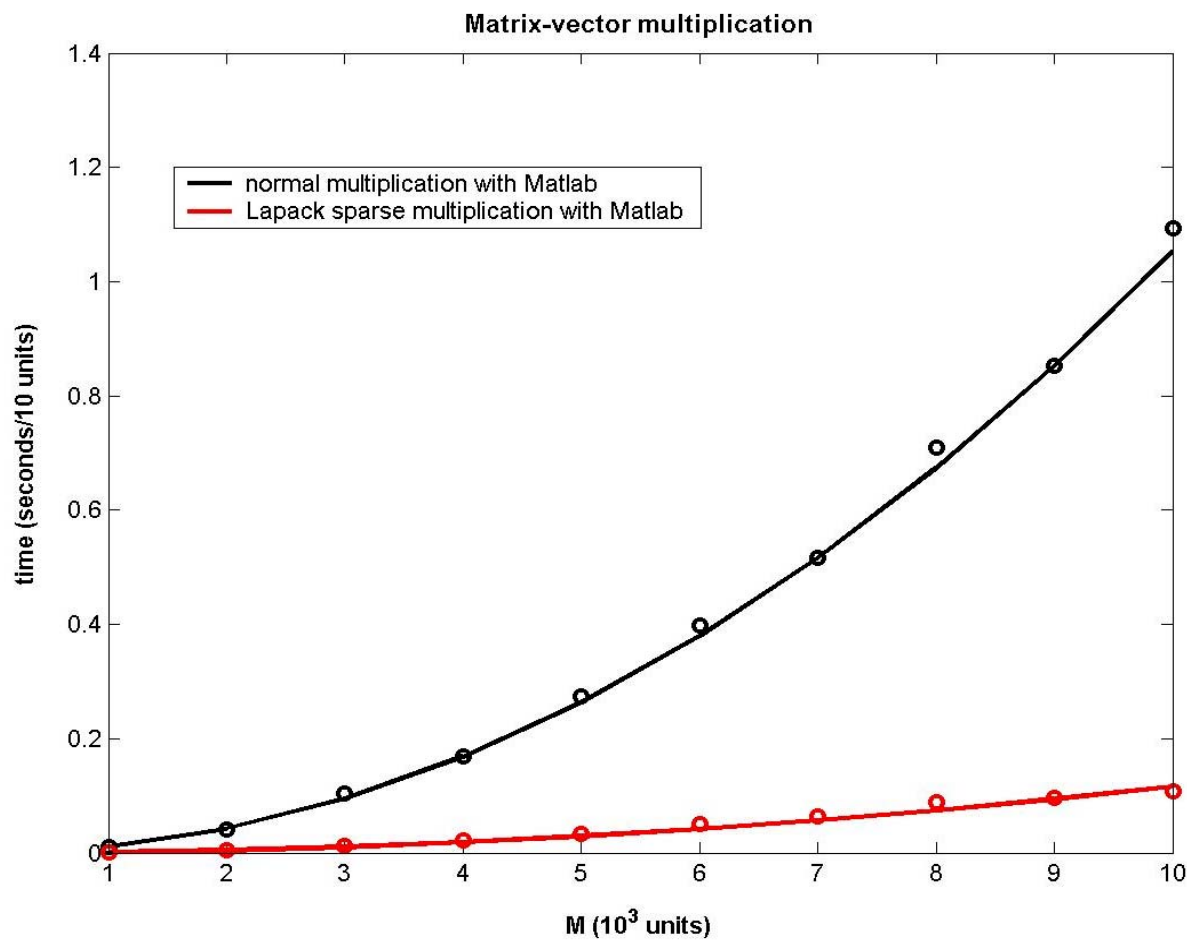

These are the results for single multiplication using Intel Xeon 3.2 GHz with 1 MB internal *cache* :

for **M**=$10^4$ the memory occupied by the sparse version of **G** is only $O(10^2)$ KB instead of theoric $O(10^6)$ KB: **G** can be stored in the processor cache

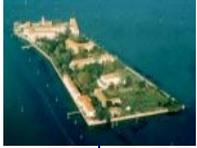

# *Distributed computing*

If we have **p** processors, with Mod(**M**,**p**)=0, we can run faster the computation of splines distributing **M/p** rows of matrix **G** to every processor:

• Single Program Multiple Data method

• no communication among processes is involved; there is only a limited overhead for sending the rows of **G** to every processor

• tests with a Matlab multi-engine environment

From previous example, using 4 Xeon processors and **N**=$10^3$, the registered execution time for computing all the splines is about 12 seconds (0.012 * 3 * $10^3$ / 4 = 9 the theoric):

$$\text{time}_{execution} = \text{time}_{cpu} + \text{time}_{overhead}$$

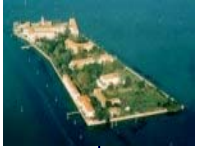

# Post-processing phase

Now we would a fast method for computing the splines values in a set of *parameter ticks* with fine sampling

Let **V** + 1 the number of ticks for each cubic spline valuation; then the ticks are (0, 1/ **V**, 2/ **V**, . . ., (**V** -1)/ **V** , 1), and the values of parameter in the computation are their (0, 1, 2, 3)-th degree powers. The value of a cubic at $t_0$ can be view as a dot product:

$at_0^3 + bt_0^2 + ct_0 + d = (a, b, c, d) \bullet (t_0^3, t_0^2, t_0, 1)$

This fact permits to consider a new *constant* 4 x (**V**+1) **T** matrix:

$$T = \begin{pmatrix} 0 & (1/\mathbf{V})^3 & . & . & ((\mathbf{V}-1)/\mathbf{V})^3 & 1 \\ 0 & (1/\mathbf{V})^2 & . & . & ((\mathbf{V}-1)/\mathbf{V})^2 & 1 \\ 0 & (1/\mathbf{V})^1 & . & . & ((\mathbf{V}-1)/\mathbf{V})^1 & 1 \\ 1 & 1 & . & . & 1 & 1 \end{pmatrix}$$

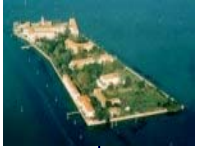

# An eulerian view

$$C = \begin{pmatrix} a_1 & b_1 & c_1 & d_1 \\ a_2 & b_2 & c_2 & d_2 \\ . & . & . & . \\ . & . & . & . \\ a_M & b_M & c_M & d_M \end{pmatrix}$$

This is the **M** x 4 matrix **C**, each row is a spline between two points, and this for all the **M** trajectories

Then the **M** x (**V**+1) matrix product **V** = **C T** contains in each row the values of a cubic between two data-points, for all the **M** trajectories (*eulerian* method: computation of all the trajectories at a predefined set of time ticks).

The theoric number of floating point operations for computing all the cubics values for all the trajectories in flow is *O(10* **MVN***)*.

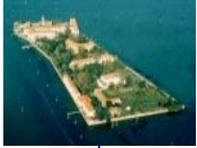

# Computing the values of splines

The matrices product is fast with Matlab incorporated Lapack routines: tests with Xeon 3.2 GHz processor, **M**=$10^4$ and **V**=10 show a time of 0.02 seconds for one multiplication;

if **N**=$10^2$, the time for computing the values of all the splines of a FD time-step (a snapshot of flow) is 0.02 x 3 x $10^2$ = 6 seconds

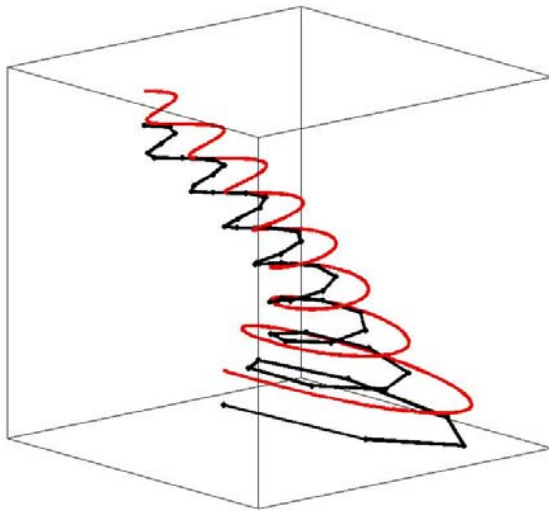

The *black line* is a real trajectory from FD computation; the *red* is the virtual line from splines method

(at left the line has been shifted; trajectory of a gas particle in combustion chamber exiting from forced ventilation fan)

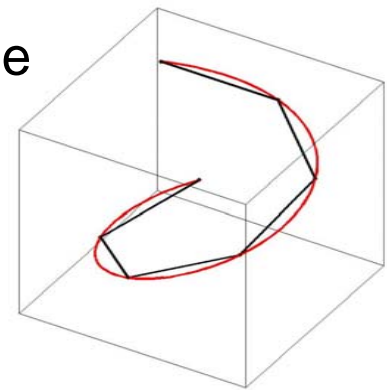

# A virtual equivalence

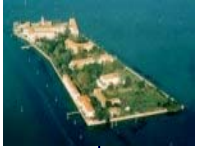

*Assume* that the splines method is *equivalent* to a Finite Difference method with a grid space-step defined by the value of **V**; from CFL:

time-step$_{splines}$ = **L** x (3**N**)$^{-1}$ x **s**$^{-1}$ , where **L** is the linear length and **s** the scalar speed of flow

time-step$_{FD}$ = **L** x (3**NV**)$^{-1}$ x **s**$^{-1}$

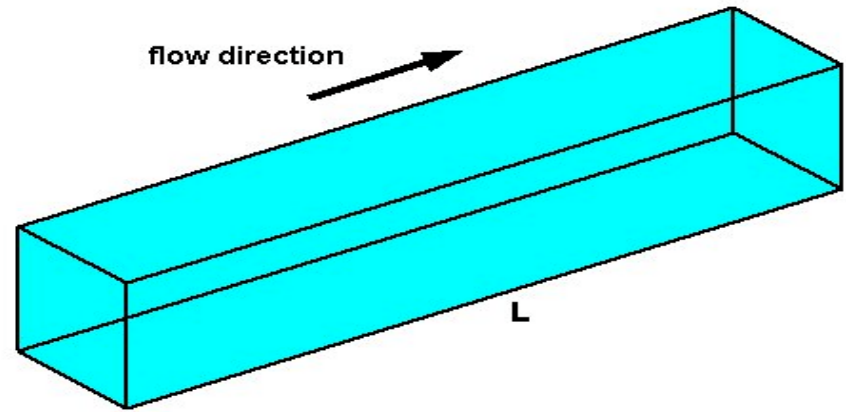

**Example** with **L**=3m, **M**=10$^4$, **N**=10$^2$, **s**=30cm/sec, **V**=10;

Report of total execution time (FD computation & graphics) for 1 minute of real simulation, 1 second between two snapshots, with four 3.2 GHz Xeon processors, for the two methods (Matlab parallel multi-engine environment):

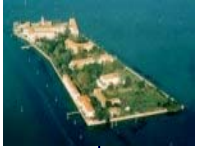

# Numerical considerations

$time_{splines}$ ~ 450 sec (FDc) + 4500 sec (splines) + 3000 sec (graphic rendering) ~ 8000 sec

$RAM_{splines}$ (total allocation) ~ 1 GB for every flow snapshot

$time_{FD}$ ~ 8000 sec (FDc) + 2000 sec (graphic rendering) ~ 10000 sec

$RAM_{FD}$ (total allocation) ~ 2 GB for every flow snapshot

**Example** of graphic output for a **V**=4 grid step (~ 2 mm space-step)

red = line from splines method (AB first spline, BC second spline)

black = points from full FD method

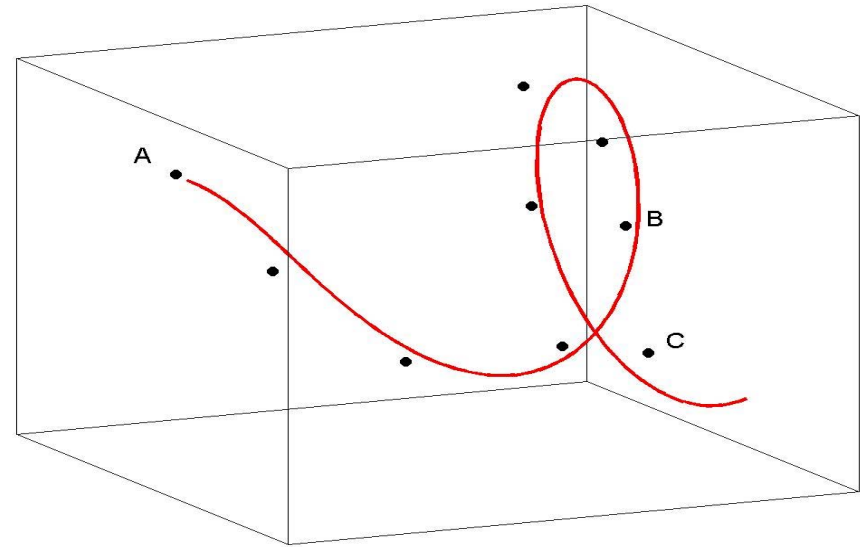

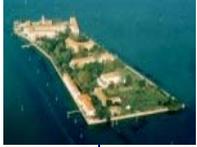

# *Pros and cons*

The **spline method** :

• reduces total time of computation and RAM allocation

• is easily adaptable for a multiprocessor architecture

• its graphic rendering gives results comparable with those of FD method

• *it is limited to the design of particles trajectories and gives no information on other Navier-Stokes variables as pressure, temperature*

• *it is not useful in the case of very small geometries (the interpoled trajectories can cut the small element of boundary)*

• *its equivalence with a finer-grid FD method must to be mathematically justified*

Thank you

September, 2004